\documentclass[12pt,a4wide]{article}
\usepackage{epsfig}
\usepackage{amsmath,amssymb,amsthm}
\usepackage{hyperref}

\theoremstyle{plain}
\newtheorem{lemma}{Lemma}[section]
\newtheorem{theorem}[lemma]{Theorem}
\newtheorem{prop}[lemma]{Proposition}
\newtheorem{cor}[lemma]{Corollary}
\newtheorem{claim}{Claim}
\theoremstyle{definition}
\newtheorem{remark}[lemma]{Remark}
\newtheorem{Ex}[lemma]{Example}
\newtheorem*{Q}{Question}

\DeclareMathOperator\supp{supp}
\DeclareMathOperator\wt{wt}
\DeclareMathOperator\Ev{Spec}
\DeclareMathOperator\nullsp{nullsp}

\newcommand{\bma}{\mathbb{A}}

\newcommand{\zero}{\mathbf{0}}
\newcommand{\one}{\mathbf{1}}

\newcommand{\vd}{\mathbf{d}}
\newcommand{\ve}{\mathbf{e}}
\newcommand{\x}{\mathbf{x}}
\renewcommand{\u}{\mathbf{u}}
\newcommand{\vv}{\mathbf{v}}

\begin{document}

\title{Arithmetic completely regular codes}

\author{J.~H.~Koolen\footnote{
School of Mathematical Sciences, University of Science and Technology, and Wu Wen-Tsun Key Laboratory of Mathematics of CAS, 
Hefei Anhui 230026, P.R.~China} 
\and W.~S.~Lee\footnote{Department of Mathematics, POSTECH, Pohang, South
Korea}
 \and W.~J.~Martin\addtocounter{footnote}{0}\footnote{
Department of Mathematical Sciences and Department of Computer
Science, Worcester Polytechnic Institute, Worcester,
Massachusetts, USA} 
\and H.~Tanaka\addtocounter{footnote}{0}\footnote{Research Center for Pure and Applied Mathematics,
Graduate School of Information Sciences,
Tohoku University, Sendai 980-8579, Japan}
}

\date{}

\maketitle

\begin{abstract}
In this paper, we explore completely regular codes in the Hamming
graphs and related graphs.  Experimental evidence suggests that many completely regular
codes have the property that the eigenvalues of the code are in arithmetic progression. 
In order to better understand these ``arithmetic completely regular codes'', we focus 
on cartesian products of completely regular codes and products of their corresponding 
coset graphs in the additive case. Employing earlier results, we are then able to prove 
a theorem which nearly classifies these codes in the case where the graph admits a 
completely regular partition into such codes (e.g, the cosets of some additive 
completely regular code).  Connections to the theory of distance-regular graphs are 
explored and several open questions are posed.
\end{abstract}

\section{Introduction}
\label{Sec:intro}

In this paper, we present the theory of completely regular codes
in the Hamming graph enjoying the property that the eigenvalues of
the code are in arithmetic progression. We call these codes {\em
arithmetic completely regular codes} and we classify them under some
additional conditions.  Our results are strongest when the Hamming graph
admits a completely regular partition into such codes (e.g., the
partition into cosets of some additive completely regular code),
since it is known that the quotient graph obtained from any such
partition is distance-regular. Using Leonard's Theorem, the list
of possible quotients is determined, with a few special cases left
as open problems. In the case of linear arithmetic completely
regular codes, more can be said.

Aside from the application of Leonard's Theorem, the techniques 
employed are mainly combinatorial and products of codes as well 
as decompositions of codes into ``reduced'' codes play a fundamental role.
Indeed, whenever a product of two completely regular codes is completely 
regular, all three codes are necessarily arithmetic and taking cartesian products
with entire Hamming graphs does not affect either the completely regular or 
the arithmetic property.

These results are not so relevant to classical coding theory;
they belong more to the theory of distance-regular graphs. 
Coset graphs of additive completely regular codes provide an
important family of distance-regular graphs (and dual pairs of 
polynomial association schemes). While the local structure of an
arbitrary distance-regular graph is hard to recover from its
spectrum alone, these coset graphs naturally inherit much of their 
combinatorial structure from the Hamming graphs from which they arise.
What is interesting here is how this additional information gives 
leverage in the combinatorial analysis of these particular distance-regular
coset graphs, a tool not available in the unrestricted case. Recent related work appears in \cite{RZ}.

It is important here to draw a connection between the present
paper and a companion paper \cite{KLM}.
In \cite{KLM}, the first three authors introduced the concept of {\em Leonard
completely regular codes} and
developed the basic theory for
them. We also gave several important families of Leonard
completely regular codes in the classical distance-regular graphs,
demonstrating their fundamental role in structural questions about
these graphs. In fact, we wonder if there exist any completely regular 
codes in the Hamming graphs with large covering radius that are neither
Leonard nor closely related to Leonard codes. The class of arithmetic 
codes we introduce in this paper
is perhaps the most important subclass of the Leonard completely
regular codes in the Hamming graphs and something similar is
likely true for the other classical families, but this
investigation is left as an open problem.

The layout of the paper is as follows. After an introductory
section outlining the required background, we explore products of
completely regular codes in Hamming graphs.  First noting
(Proposition~\ref{PCxC'crimplesbothcr}) that a completely regular
product must arise from completely regular constituents, we
determine in Proposition \ref{Pprodcrc} exactly when the product
of two completely regular codes in two Hamming graphs is
completely regular. At this point, the role of the arithmetic
property becomes clear and we understand that  
Lemma \ref{Ltridiageigenvalues} gives a generic form for the
quotient matrix of such a code.

With this preparatory material out of the way, we are then ready
to present the main results in Section
\ref{Subsec:classification}. Applying the celebrated theorem of
Leonard, Theorem \ref{Tarithmeticquo} imposes powerful
limitations on the combinatorial structure of a  quotient of a
Hamming graph when the underlying completely regular partition is
composed of arithmetic codes. Stronger results are obtained in
Proposition \ref{Pcliques}, Theorem \ref{Tbreakdown} and
Proposition \ref{PnotDoob} when one makes additional assumptions
about the minimum distance of the codes or the specific structure
of the quotient. When $C$ is a linear completely regular code with
the arithmetic property and $C$ has minimum distance at least
three and covering radius at most two, we show that $C$ is closely
related to some Hamming code. These results are summarized in
Theorem \ref{Tproductrho=1}, which gives a full classification of
possible codes and quotients in the linear case (always assuming
the arithmetic and completely regular properties) and Corollary
\ref{CHammingquo} which characterizes Hamming quotients of Hamming
graphs.

\section{Preliminaries and definitions}
\label{Sec:prelim}

In this section, we summarize the background material necessary to
understand our results. Most of what is covered here is based on
Chapter 11 in the monograph \cite{B} by Brouwer, Cohen and
Neumaier. The theory of codes in distance-regular graphs began
with Delsarte \cite{D}. The theory of association schemes is also
introduced in the book \cite{A} of Bannai and Ito while
connections between these and related material (especially
equitable partitions) can be found in Godsil \cite{Godsil}.

\subsection{Distance-regular graphs}

Suppose that $\Gamma$ is a finite, undirected, connected graph
with vertex set $V\Gamma$. For vertices $x$ and $y$ in $V\Gamma$,
let $d(x,y)$ denote the distance between $x$ and $y$, i.e., the
length of a shortest path connecting $x$ and $y$ in $\Gamma$. Let
$D$ denote the diameter of $\Gamma$; i.e., the maximal distance
between any two vertices in $V\Gamma$. For $0\le i\le D$ and $x\in
V\Gamma$, let ${\Gamma_{i}}(x):=\{y\in V \Gamma \,|\, d(x,y)=i\}$
and put $\Gamma_{-1}(x):=\emptyset$, $\Gamma_{D+1}(x):=\emptyset$.
The graph $\Gamma$ is called {\em distance-regular} whenever it is
regular of valency $k$, and there are integers $b_{i},c_{i} \
(0\leq i \leq D)$ so that for any two vertices $x$ and $y$ in
$V\Gamma$ at distance $i$, there are precisely $c_{i}$ neighbors
of $y$ in $\Gamma_{i-1}(x)$ and $b_{i}$ neighbors of $y$ in
$\Gamma_{i+1}(x)$.  It follows that there are exactly
$a_i=k-b_i-c_i$ neighbors of $y$ in $\Gamma_i(x)$. The numbers
$c_{i}$, $b_{i}$ and $a_{i}$ are called the {\em intersection
numbers} of $\Gamma$ and we observe that $c_{0}=0$, $b_{D}=0$,
$a_{0}=0$, $c_{1}=1$ and $b_{0}=k$. The array
$\iota(\Gamma):=\{b_{0},b_{1},\dots,b_{D-1};c_{1},c_{2},\dots,c_{D}\}$
is called the {\em intersection array} of $\Gamma$.

From now on, assume $\Gamma$ is a distance-regular graph of
valency $k\ge 2$ and diameter $D\ge 2$. Define $A_{i}$ to be the
square matrix of size $|V\Gamma|$  whose rows and columns are
indexed by $V\Gamma$ with entries
\begin{equation*}
	(A_{i})_{xy}= \begin{cases} 1 & \text{if} \ d(x,y)=i \\
0 & \text{otherwise} \end{cases} \quad (0\leq i\leq
D, \ x,y \in V\Gamma).
\end{equation*}

We refer to $A_{i}$ as the {\em $i^{\rm th}$ distance matrix} of
$\Gamma$. We abbreviate $A:=A_{1}$ and call this the {\em
adjacency matrix} of $\Gamma$. Since $\Gamma$ is distance-regular,
we have
\begin{equation*}
	A A_{i-1} = b_{i-2} A_{i-2} + a_{i-1} A_{i-1} + c_i A_i \quad (2\le i \le D)
\end{equation*}
so that $A_{i}=p_{i}(A)$ for some polynomial $p_{i}(t)$ of degree
$i$. Let $\bma$ be the {\em Bose-Mesner algebra}, the matrix
algebra over $\mathbb{C}$ generated by $A$. Then $\dim \bma =D+1$
and $\{A_{i}\,|\,0 \le i \le D\}$ is a basis
for $\bma$. As $\bma$ is semi-simple and commutative, $\bma$ has
also a basis of pairwise orthogonal idempotents $\left\{
E_{0}=\frac{1}{|V\Gamma|}J,E_{1},\dots,E_{D} \right\}$, where $J$ denotes the all ones matrix.
We call these matrices the {\em primitive idempotents} of $\Gamma$. As
$\bma$ is closed under the entrywise (or Hadamard or Schur) product
$\circ$, there exist real numbers $q_{ij}^\ell$, called the {\em
Krein parameters}, such that
\begin{equation}
\label{Ekrein} E_i\circ E_j =\frac{1}{|V\Gamma|}\sum_{\ell=0}^{D}
q_{ij}^\ell E_\ell \quad ( 0\leq i,j \leq D).
\end{equation}

The graph $\Gamma$ is called {\em $Q$-polynomial} if there exists
an ordering $E_{0},\dots,E_{D}$ of the primitive idempotents and
there exist polynomials $q_{i}$ of degree $i$ such that
$E_{i}=q_{i}(E_{1})$, where the polynomial $q_i$ is applied
entrywise to $E_1$. We recall that the distance-regular graph
$\Gamma$ is $Q$-polynomial with respect to the ordering
$E_0,E_1,\dots, E_D$ provided its Krein parameters satisfy
\begin{itemize}
\item $q_{ij}^\ell =0$ unless $|j-i| \le \ell \le i+j$;
\item $q_{ij}^\ell \neq 0$ whenever $\ell=|j-i|$ or $\ell= i+j \le D$.
\end{itemize}

By an {\em eigenvalue} of $\Gamma$, we mean an eigenvalue of
$A=A_1$. Since $\Gamma$ has diameter $D$, it has at least $D+1$ distinct
eigenvalues; but since $\Gamma$ is distance-regular, it has
exactly $D+1$ distinct eigenvalues and $E_0,\dots,E_D$ are the orthogonal projections onto the eigenspaces.
We denote by $\theta_i$ the eigenvalue associated with $E_i$ and, aside from the convention that
$\theta_0=k$, the valency of $\Gamma$, we make no further
assumptions at this point about the eigenvalues except that they
are distinct.

\subsection{Codes in distance-regular graphs}

Let $\Gamma$ be a distance-regular graph with distinct eigenvalues
$\theta_{0}=k,\theta_{1},\dots,\theta_{D}$. By a {\em code} in
$\Gamma$, we simply mean any nonempty subset $C$ of $V \Gamma$. We
call $C$ {\em trivial} if $|C| =1$ or $C=V\Gamma$ and {\em
non-trivial} otherwise. For $|C|>1$, the {\em minimum distance} of
$C$, $\delta(C)$, is defined as
\begin{equation*}
	\delta(C):= \min \{\,d(x,y) \,|\, x,y\in C,x\neq y\,\}
\end{equation*}
and for any $x \in V \Gamma $ the distance $d(x,C)$ from $x$ to 
$C$ is defined as
\begin{equation*}
	d(x,C):= \min \{\,d(x,y) \,|\, y\in C\,\}.
\end{equation*}
The number
\begin{equation*}
	\rho(C):= \max \{\,d(x,C) \,|\, x\in V \Gamma\,\}
\end{equation*}
is called the {\em covering radius} of $C$.

For a code $C$
and for $0\leq i \leq\rho:=\rho(C)$,
define
\begin{equation*}
	C_{i}=\{\,x\in V\Gamma \,|\, d(x,C)=i\,\}.
\end{equation*}
Then $\Pi(C):=\{C_{0}=C,C_{1},\dots,C_{\rho}\}$ is the {\em
distance partition} of $V \Gamma$ with respect to code $C$.

A partition $\Pi=\{P_{0},P_{1},\dots,P_{\ell}\}$ of $V\Gamma$ is
called {\em equitable} if, for all $i$ and $j$, the number of
neighbors a vertex in $P_{i}$ has in $P_{j}$ is independent of the
choice of vertex in $P_{i}$.  Following Neumaier~\cite{G}, we say
a code $C$ in $\Gamma$ is {\em completely regular} if this
distance partition $\Pi(C)$ is equitable\footnote{This definition
of a completely regular code is equivalent to the original
definition, due to Delsarte \cite{D}.}. In this case the following
quantities are well-defined:
\begin{align}
	\gamma_{i} &= \left|\{y\in C_{i-1}\,|\,d(x,y)=1\} \right|, \label{Egamma} \\
	\alpha_{i} &= \left|\{y\in C_{i}\,|\,d(x,y)=1\}\right|, \label{Ealpha} \\
	\beta_{i} &= \left|\{y\in C_{i+1}\,|\,d(x,y)=1\}\right| \label{Ebeta}
\end{align}
where $x$ is chosen from $C_{i}$. The numbers
$\gamma_{i},\alpha_{i},\beta_{i}$ are called the {\em intersection
numbers} of code $C$. Observe that a graph $\Gamma$ is
distance-regular if and only if each vertex is a completely
regular code and these $|V\Gamma|$ codes all have the same
intersection numbers. Set the tridiagonal matrix
\begin{equation*}
	U=U(C):=\begin{pmatrix}
	\alpha_{0} & \beta_{0}  &    &         &  \\
	\gamma_{1} & \alpha_{1} & \beta_{1} &  &   \\
	   & \gamma_{2} & \alpha_{2} &\beta_{2}  &   \\
	   &            &  \ddots & \ddots  &  \ddots  \\
	   &           & & \gamma_{\rho}&\alpha_{\rho}\end{pmatrix}.
\end{equation*}

For $C$ a completely regular code in $\Gamma$, we say that $\eta$
is an {\em eigenvalue of $C$} if $\eta$ is an eigenvalue of the
{\em quotient matrix} $U$ defined above. By $\Ev(C)$, we denote the set
of eigenvalues of $C$. Note that, since
$\gamma_i+\alpha_i+\beta_i=k$ for all $i$, $\theta_0=k$ belongs to
$\Ev(C)$.

\subsection{Completely regular partitions}

Given a partition $\Pi$ of the vertex set of a graph $\Gamma$
(into nonempty sets), we define the {\em quotient graph}
$\Gamma/\Pi$ on the cells of $\Pi$ by calling two cells
$C,C'\in\Pi$ adjacent if $C\neq C'$ and $\Gamma$ contains an edge
joining some vertex of $C$ to some vertex of $C'$. A partition
$\Pi$ of $V\Gamma$ is {\em completely regular} if it is an
equitable partition of $\Gamma$ and all $C\in \Pi$ are completely
regular codes with the same intersection numbers.
If $\Pi$ is completely regular we write $\Ev(\Pi)=\Ev(C)$ for $C\in\Pi$ and we say $\rho:=\rho(C)$ is the covering radius of $\Pi$.

\begin{prop}[cf.~{\cite[p352-3]{B}}]
\label{Pbcnquo}
Let $\Pi$ be a completely regular partition of any distance-regular graph
$\Gamma$ such that each $C\in\Pi$ has intersection numbers
$\gamma_i$, $\alpha_i$ and $\beta_i$ $(0\le i \le \rho)$. Then
$\Gamma/\Pi$ is a distance-regular graph with intersection array
\begin{equation*}
	\iota(\Gamma/\Pi) = \left\{
	\frac{\beta_{0}}{\gamma_{1}} , \frac{\beta_{1}}{\gamma_{1}} , \dots,
	\frac{\beta_{\rho-1}}{\gamma_{1}} ;\, 1,
	\frac{\gamma_{2}}{\gamma_{1}} , \dots,
	\frac{\gamma_{\rho}}{\gamma_{1}}  \right\},
\end{equation*}
remaining intersection numbers 
$a_i = \frac{\alpha_i-\alpha_0}{\gamma_1}$, and eigenvalues
 $\frac{\theta_j - \alpha_0}{\gamma_1}$ for $\theta_{j} \in {\Ev}(C)$.
All of these lie among the eigenvalues of the matrix
$ \frac{1}{\gamma_1} ( A - \alpha_0 I )$.
\end{prop}

\begin{prop}
\label{Psmallesteig} Let $\Pi$ be a non-trivial completely regular
partition of a distance-regular graph $\Gamma$ and assume that
$\Ev(\Pi)=\{\eta_{0}\geq\eta_{1}\geq\dots\geq\eta_{\rho}\}$. Then
$\eta_{\rho}\leq\alpha_{0}-\gamma_{1}$.
\end{prop}

\begin{proof}
By Proposition \ref{Pbcnquo}, the eigenvalues of $\Gamma/\Pi$ are
$\frac{\eta_{0}-\alpha_{0}}{\gamma_{1}},\frac{\eta_{1}-\alpha_{0}}{\gamma_{1}},\dots,\frac{\eta_{\rho}-\alpha_{0}}{\gamma_{1}}$.
As $\Gamma/\Pi$ has at least one edge, it follows that its
smallest eigenvalue is at most $-1$. Hence
$\frac{\eta_{\rho}-\alpha_{0}}{\gamma_{1}}\leq -1$.
\end{proof}

\begin{Q}
Let $C$ be a non-trivial completely regular
code in a distance-regular graph $\Gamma$. Let
$\theta:=\min\{\eta\,|\,\eta \in \Ev(C)\}$. Is
it true that $\theta\leq\alpha_{0}-\gamma_{1}$?
\end{Q}

\section{Codes in the Hamming graph}
\label{Sec:Hnq}

Let $X$ be a finite abelian group. A {\em translation
distance-regular graph} on $X$ is a distance-regular graph
$\Gamma$ with vertex set $X$ such that if $x$ and $y$ are adjacent
then $x+z$ and $y+z$ are adjacent for all $x,y,z\in X$. A code
$C\subseteq X$ is called {\em additive} if for all $x,y\in C$,
also $x-y \in C$; i.e., $C$ is a subgroup of $X$. If $C$ is an
additive code in a translation distance-regular graph $\Gamma$ on $X$, then
we obtain the usual coset partition $\Delta(C):=\{C+x
\,|\, x\in X \}$ of $X$;  whenever $C$ is a
completely regular code, it is easy to see that $\Delta(C)$ is a
completely regular partition.
The quotient graph $\Gamma/\Delta(C)$ is usually referred to as the \emph{coset graph} of $C$ in $\Gamma$.
An important special class of additive codes are the linear codes.
Here, assume $X$ is a vector space over some finite field $GF(q)$; a 
code $C\subseteq X$ is a {\em linear code} if $C$ is a vector subspace of $X$.
In our results, we are careful to assume only additivity
when possible, but some results only apply in the linear case.

An important result of Brouwer, Cohen and Neumaier \cite[p353]{B}
states that every translation distance-regular graph  of diameter
at least three defined on an elementary abelian group $X$ is
necessarily a coset graph of some additive completely regular code
in some Hamming graph. Of course, the Hamming graph itself is a
translation graph, usually in a variety of ways.

Let $\mathcal{Q}$ be an abelian group with $|\mathcal{Q}|=q$.  Then we may identify
the vertex set of the Hamming graph $H(n,q)$ with the group
$X=\mathcal{Q}^n$, two vertices being adjacent if they agree in all but one coordinate position.
It is well known (e.g., \cite[Theorem 9.2.1]{B})
that $H(n,q)$ has eigenvalues $\theta_j=n(q-1)-qj$, $0\le j\le n$. The corresponding 
``natural'' ordering $E_0,E_1,\dots,E_n$ is  a $Q$-polynomial ordering. When $q > 2$ and
$n>2$, this is the unique $Q$-polynomial ordering for $H(n,q)$ but when 
$q=2$ and $n$ is even,
$H(n,2)$ is also $Q$-polynomial with respect to the ordering $E_0,E_{n-1},E_2,E_{n-3},\dots, E_1,E_n$ \cite[p305]{A}.
(This second ordering plays a role in the addendum.)

For the remainder of
this paper, we will consider $q$-ary codes of length $n$ as
subsets of the vertex set of the Hamming graph $H(n,q)$. In this
section we will focus on $q$-ary completely regular codes, i.e.,
$q$-ary codes of length $n$ which are completely regular in
$H(n,q)$.

\subsection{Products of completely regular codes}
\label{Subsec:products}

We now recall the cartesian product of two graphs $\Gamma$ and
$\Sigma$. This new graph $\Gamma\times\Sigma$ has vertex set $V\Gamma \times V\Sigma$
and adjacency $\sim$ defined by $(x,y)\sim (u,v)$ precisely when
either $x=u$ and $y\sim v$ in $\Sigma$ or $x\sim u$ in $\Gamma$
and $y=v$. Now let $C$ be a code in $\Gamma$ and let $C'$ be a code in $\Sigma$.
The cartesian product of $C$ and $C'$ is the code defined by
\begin{equation*}
	C\times C' :=\left\{(x,x') \in V\Gamma \times V\Sigma \, | \, x \in C, \, x'\in C' \right\}.
\end{equation*}
We are interested in the cartesian product of codes in the Hamming 
graphs $H(n,q)$ and $H(n',q')$. Note that if $C$ and $C'$ are
additive codes then $C\times C'$ is also additive. Also,
if $C$ is a (not necessarily additive) completely regular code in $H(n,q)$ with vertex set
$\mathcal{Q}^n$, then $\mathcal{Q}\times C$ is a completely regular code in $H(n+1,q)$.
Also if $\Pi$ is a completely regular partition
of $H(n,q)$ then
$\mathcal{Q}\times \Pi:=\{\mathcal{Q}\times P\,|\,P\in\Pi\}$
is a completely regular partition of $H(n+1,q)$. We say a completely
regular code $C$ in $H(n,q)$ is {\em non-reduced} if
$C\cong \mathcal{Q}\times C'$ for some $C'$ in $H(n-1,q)$, and {\em
reduced} otherwise. In similar fashion we say that a completely
regular partition is {\em non-reduced} or {\em reduced}.

The next three results will determine exactly when the cartesian
product of two arbitrary codes in two Hamming graphs is completely
regular.

\begin{prop}
\label{PCxC'crimplesbothcr} Let $C$ and $C'$ be non-trivial codes
in the Hamming graphs $H(n,q)$ and $H(n',q')$ respectively. If the
cartesian product $C\times C'$ is completely regular in $H(n,q)
\times H(n',q')$, then both $C$ and $C'$ themselves must be
completely regular codes in their respective graphs.
\end{prop}

\begin{proof}
Let $\{C_0=C,C_1,\dots, C_\rho\}$ be the distance partition of $H(n,q)$ with respect to $C$.
Then we easily see that every
vertex $(x,y)$ of $C_i \times C'$ is at distance $i$ from $C\times
C'$ in the product graph. Moreover, the neighbors of this vertex
which lie at distance $i-1$ from the product code are precisely
$\{ (u,y) \,|\, u \sim x, \, u\in C_{i-1} \}$; so the size of the
set $\{ u \in C_{i-1} \,|\, u \sim x \}$ must be
independent of the choice of $x \in C_i$. This shows that the
intersection numbers $\gamma_i$ are well-defined for $C$. An
almost identical argument gives us the intersection numbers
$\beta_i$. This shows that $C$ is a completely regular code
and, swapping the roles of $C$ and $C'$, we show the same is true
for $C'$.
\end{proof}

It is well known that equitable partitions are preserved under
products. If $\Pi$ is an equitable partition in any graph $\Gamma$
and $\Delta$ is an equitable partition in another graph $\Sigma$,
then $\left\{ P\times P' \, | \ P \in \Pi, \ P' \in \Delta
\right\}$ is an equitable partition in the cartesian product graph
$\Gamma\times \Sigma$. The following special case will prove
useful in our next proposition.

\begin{lemma}
\label{Lequitprod}
Let $C$ be a completely regular code in $H(n,q)$ with distance partition $\Pi = \{C_0=C,\dots,
C_\rho\}$, and let $C'$ be a completely regular code in $H(n',q')$ with distance partition $\Pi' = \{C'_0=C',\dots, C'_{\rho'}\}$.
\begin{enumerate}
\item The partition
\begin{equation*}
	\left\{ C_i \times C'_j \,|\, 0\le i\le \rho, \, 0\le j\le \rho'\right\}
\end{equation*}
is an equitable partition in the product  graph $H(n,q) \times
H(n',q')$.
\item A vertex in $C_i \times C'_j$ has all its neighbors in
\begin{equation*}
	\left[ \left( C_{i-1} \cup C_i \cup C_{i+1} \right) \times C'_j \right]
	\cup \left[ C_i \times \left( C'_{j-1} \cup C'_j \cup C'_{j+1} \right)\right].
\end{equation*}
\item If $C$ has intersection numbers $\gamma_i,\alpha_i,\beta_i$
and $C'$ has intersection numbers $\gamma'_i,\alpha'_i,\beta'_i$, then
in the product graph, a vertex in $C_i \times C'_j$ has:
$\gamma_i$ neighbors in $C_{i-1} \times C'_j$; $\beta_i$ neighbors
in $C_{i+1} \times C'_j$; $\gamma'_j$ neighbors in $C_i \times
C'_{j-1}$; $\beta'_j$ neighbors in $C_i \times C'_{j+1}$; and
$\alpha_i + \alpha'_j$ neighbors in $C_i \times C'_j$.
\end{enumerate}
\end{lemma}

\begin{proof}
Straightforward.
\end{proof}

\begin{lemma}
\label{Ltridiageigenvalues} Let $k,\gamma,\beta$ and $\rho$ be
positive integers. The tridiagonal matrix
\begin{equation*}
	L=\begin{pmatrix}
	\alpha_{0} & \rho\beta  &               &              &           &  \\
	    \gamma & \alpha_{1} & (\rho-1)\beta &              &           &  \\
	            & 2\gamma   & \alpha_{2}    &(\rho-2)\beta &           &  \\
	            &           & \ddots        & \ddots       & \ddots    &  \\
	            &           &               & (\rho-1)\gamma& \alpha_{\rho-1} & \beta \\
	            &           &               &              &\rho\gamma &\alpha_{\rho}
	\end{pmatrix}
\end{equation*}
where $\alpha_{i}=k-i\gamma-(\rho-i)\beta$ $(0\leq
i\leq\rho)$, has eigenvalues
$\Ev(L)=\{k-ti\,|\,0\leq i\leq\rho\}$ where
$t=\gamma+\beta$.
\end{lemma}

\begin{proof}
By direct verification. (This is not new. See, for example, type (IIC) in Terwilliger \cite[Theorem 2.1]{T}.)
\end{proof}

\begin{prop}
\label{Pprodcrc} Let $C$ be a non-trivial completely regular code
in $H(n,q)$  with $\rho:=\rho(C)\geq 1$ and intersection numbers
$\alpha_i,\beta_i$ and $\gamma_i$ $(0\le i \le \rho)$. Let $C'$ be
a non-trivial completely regular code in $H(n',q')$  with
$\rho':=\rho(C')\geq 1$ and intersection numbers
$\alpha'_i,\beta'_i$ and $\gamma'_i$ $(0\le i \le \rho')$.
Then $C\times C'$ is a
completely regular code in $H(n,q) \times H(n',q')$ if and only if
there exist integers $\gamma$ and $\beta$ satisfying (a) and (b):
\begin{enumerate}
\item $\gamma_i = \gamma i$ for $0\le i \le \rho$ and $\gamma'_i = \gamma i$ for $0\le i \le \rho'$;
\item $\beta_{\rho-i} = \beta i$ for $0\le i \le \rho$ and $\beta'_{\rho'-i} = \beta i$ for $0\le i \le \rho'$.
\end{enumerate}
In this case, $C\times C'$ has covering radius $\bar{\rho}:=
\rho+\rho'$ and intersection numbers $\bar{\gamma}_i = \gamma i$ and
$\bar{\beta}_i = \beta(\bar{\rho}-i)$ for $0\le i  \le \bar{\rho}$, and all
three codes --- $C$, $C'$ and $C\times C'$ --- are arithmetic completely regular codes.
\end{prop}

\begin{proof}
Assume, without loss, that $\rho \le \rho'$. 
$(\Rightarrow)$ Suppose first that $C\times C'$ is a completely
regular code. From Lemma \ref{Lequitprod}, we see that $C\times
C'$ has covering radius $\bar{\rho}=\rho+\rho'$. For $0\le j\le
\bar{\rho}$, let
\begin{equation*}
	S_j = \left\{  (x,y) \,|\, d( (x,y),C\times C') = j  \right\}.
\end{equation*}
Then it follows easily from Lemma \ref{Lequitprod}(b) that
\begin{equation*}
	S_j = \bigcup_{ h+i=j} C_h \times C'_i .
\end{equation*}
Moreover, by Lemma \ref{Lequitprod}(c), a vertex in $C_h \times C'_i$ has $\gamma_h + \gamma'_i$
neighbors in $S_{j-1}$ and $\beta_h + \beta'_i$ neighbors in
$S_{j+1}$.  For $j=1$, this forces $\gamma_1=\gamma'_1=:\gamma$. Assume
inductively that $\gamma_i = \gamma i$  and $\gamma'_i =
\gamma i$ for $i<j$. Then, considering a vertex in
\begin{equation*}
	S_j = C_r \times C'_{j-r} \cup \dots \cup C_{j-s} \times C'_s
\end{equation*}
(where $r = \max\{ 0, j-\rho'\}$ and $s=\max\{0,j-\rho\}$), we find
\begin{equation*}
	\gamma_r + \gamma'_{j-r} = \gamma_{r+1} + \gamma'_{j-r-1} = \dots =
	\gamma_{j-s} + \gamma'_{s}.
\end{equation*}
For $j\le \rho$, this gives $\gamma_j=\gamma j = \gamma'_j$.
For $\rho < j \le \rho'$, we deduce, $\gamma'_j=\gamma j$.
So we have (a) by induction. A symmetrical argument establishes
part (b).

$(\Leftarrow)$ Considering the same partition of $S_j$ into cells
of the form $C_i \times C'_{j-i}$, we obtain the converse in a
straightforward manner: if $C$ and $C'$ have intersection numbers
given by (a) and (b), then their cartesian product $C\times C'$ is
completely regular in the product graph. 

The fact that all three codes are arithmetic completely regular codes when (a)
and (b) hold now follows directly from Lemma \ref{Ltridiageigenvalues}.
\end{proof}

\begin{remark}
Note that, if $q\neq q'$, then the product graph is not
distance-regular. Although this case is not of primary interest,
there are examples where $C\times C'$ can still be a completely
regular code (in the sense of Neumaier) in such a graph. For
instance, suppose that $C$ is a perfect code with covering radius one in
$H(n,q)$ and $C'$ is a perfect code with covering radius one in
$H(n',q')$. If we happen to have $n(q-1)=n'(q'-1)$, then, by the
above proposition, $C\times C'$ is a completely regular code in
$H(n,q) \times H(n',q')$ with intersection numbers $\gamma_1=1$, $\gamma_2=2$, $\beta_0=2n(q-1)$ and $\beta_1=n(q-1)$.
\end{remark}

From now on we will look at a completely regular partition $\Pi$ of
$H(n,q)$ into arithmetic completely regular codes
with $\Ev(\Pi)=\{n(q-1),n(q-1)-qt, \dots ,n(q-1)-q\rho t\}$.
As a direct consequence of Proposition
\ref{Psmallesteig}, we have the following proposition which says that $\alpha_0$ must be quite large 
unless at least one of $\rho$ (the covering radius) or $t$ (the spectral gap) is large.

\begin{prop}
Let $\Pi$ be a non-trivial completely regular partition of the
Hamming graph $H(n,q)$ such that $\Pi$ has covering radius $\rho$
and $\Ev(\Pi)=\{n(q-1),n(q-1)-qt, \dots ,n(q-1)-q\rho t\}$ for
some $t$. Then $n(q-1)-\alpha_{0}\leq q\rho t$.
\end{prop}

\subsection{Classification}
\label{Subsec:classification}

In this section, we will employ  results from \cite{KLM}, in which 
{\em Leonard completely regular codes} are defined and studied. 
But we must first address an 
error in that paper.
Lemma 6.3 in \cite{KLM} claims that every harmonic completely regular code is Leonard. This fails to hold in some cases where the $Q$-polynomial ordering does not place the eigenvalues
of $\Gamma$ in decreasing order. See the addendum below for a discussion of these counterexamples
as well as a corrected statement of that lemma along with a full proof.
Since arithmetic completely regular codes are harmonic completely regular codes in the Hamming graphs (with respect to the natural ordering),
Lemma \ref{Lemma 6.3} in the addendum applies and all arithmetic completely regular codes are indeed Leonard codes.
 
Now we will classify the possible quotients $\Gamma/\Pi$ where
$\Pi$ is  a completely regular partition of the Hamming graph $\Gamma=H(n,q)$ into arithmetic
completely regular codes.
Recall that a  {\em Doob graph} \cite[p27]{B} is any
graph formed as a cartesian product of $s\ge 1$ copies of the Shrikhande
graph \cite[p104-5]{B} and $t\ge 0$ copies of $K_4$. A graph so constructed 
has the same parameters as the Hamming graph $H(2s+t,4)$ but quite
different local structure; e.g., the graph contains a four-clique if and only if $t\neq 0$ but always 
contains maximal cliques of size three. The Doob graphs are therefore important in the theory 
of distance-regular graphs \cite[Section 9.2B]{B}.

\begin{theorem}
\label{Tarithmeticquo}
Let $\Gamma$ be the Hamming graph $H(n,q)$,
and let $\Pi$ be a completely regular partition of $V\Gamma$.
Assume that $\Pi$ has covering radius $\rho\geq 3$ and eigenvalues
$\Ev(\Pi)=\{n(q-1),n(q-1)-qt, \dots ,n(q-1)-q\rho t\}$ for some
$t$. Then $\Gamma/\Pi$ has diameter $\rho$ and is isomorphic to
one of the following:
\begin{enumerate}
\item a folded cube;
\item a Hamming graph;
\item a Doob graph;
\item a distance-regular incidence graph of some $2$-$(16,6,2)$ design 
(i.e., with intersection array $\{6,5,4;1,2,6\}$).
\end{enumerate}
\end{theorem}

\begin{proof}
Denote $\eta_{i}:=k-i\tau$ $(0\leq i\leq \rho)$,  where
$k:=\frac{n(q-1)-\alpha_{0}}{\gamma_{1}}$ and
$\tau:=\frac{qt}{\gamma_{1}}$. Then by Proposition \ref{Pbcnquo},
$\Gamma/\Pi$ is a distance-regular graph with eigenvalues
$\{\eta_{i} \,|\, 0\leq i \leq\rho\}$ and
diameter $\rho$. As stated above, each code in partition $\Pi$ is a 
Leonard completely regular code. So we can apply the proof of 
\cite[Proposition~4.5]{KLM} to conclude that the quotient graph $\Gamma/\Pi$
is $Q$-polynomial with respect to $\overline{E}_{0},
\overline{E}_{1},\dots,\overline{E}_{\rho}$ where $\overline{E}_{i}$ is the primitive idempotent of $\Gamma/\Pi$
corresponding to $\eta_{i}$. (Note that the statement of Proposition~4.5 in \cite{KLM}  
assumes that $\Pi$ is a coset partition, but the proof applies in the more 
general situation we face here.)

By Leonard's Theorem \cite[p263]{A} (cf.~\cite{Leonard}),
$Q$-polynomial distance-regular graphs
with diameter $D\geq 3$ fall into seven types, namely, (I), (IA), (II), (IIA), (IIB), (IIC), and (III).
See also \cite[Section 2]{T}.
As the eigenvalues of $\Gamma/\Pi$ are in
arithmetic progression and $\rho\geq 3$, neither type
(I), (IA),  (II), (IIA) nor (III) can occur. The only remaining
possibilities are types (IIB) and (IIC).
Terwilliger \cite{I} showed that a $Q$-polynomial distance-regular
graph of type (IIB) with
 diameter $D\geq 3$ is either the folded
 $(2D+1)$-cube, or has the same intersection numbers as
 the folded $2D$-cube.
By \cite[Theorem 9.2.7]{B}, it follows that $\Gamma/\Pi$ is either
a folded $n$-cube for $n\geq 7$ or has intersection
array $\{6,5,4;1,2,6\}$.
For type (IIC), we can easily check (using $c_1=1$) that
$\Gamma/\Pi$ has the same parameters as a Hamming graph (cf.~\cite[Proposition 6.2]{Tanaka}).
But Egawa \cite{E} showed that a distance-regular graph with the same
parameters as a Hamming graph must be a Hamming graph or a Doob
graph. Hence the theorem is proved.
\end{proof}

\begin{remark}
\begin{description}
\item[A.]
Theorem \ref{Tarithmeticquo} also holds when $\Gamma$ is a Doob
graph. There are completely regular partitions of Doob graphs whose 
corresponding quotient graphs
are Hamming graphs with $q\neq 4$. For example, let $s\geq 2$ be
an integer. In any Doob graph $\Gamma$ of diameter at least five containing at least one four-clique,  there exists an additive completely regular
code, say $C$, with covering radius one and sixteen cosets \cite{F}. Let $n\geq
1$ be an integer. Then $C^n$ is an additive completely regular
code with covering radius $n$ in $\Gamma^n$,
and the coset graph $\Gamma^n/\Delta(C^n)$ is isomorphic to $H(n,16)$.
Moreover, when $\Gamma$ contains a subgraph isomorphic to $H((4^s-1)/3,4)$, we may also obtain $H(n,4^s)$ as a quotient.
\item[B.] There are three nonisomorphic graphs with the same intersection array,
$\{6,5,4;1,2,6\}$, as the folded $6$-cube. They are the
point-block incidence graphs of the $2$-$(16,6,2)$ designs. We do not
know any example in which the quotient graph has intersection
array $\{6,5,4;1,2,6\}$, but is not the folded $6$-cube.
\item[C.] If $\Gamma$ is a Hamming graph and $C$ is an additive completely regular code in
$\Gamma$ satisfying the eigenvalue conditions of Theorem \ref{Tarithmeticquo}, then $\Gamma/\Delta(C)$ is always a Hamming graph, a Doob
graph or a folded cube.  We do not know any example where $\Gamma$
is a Hamming graph and $\Pi$ is a completely regular partition
such that $\Gamma/\Pi$ is a Doob graph. But we will prove below in
Proposition \ref{PnotDoob} that this cannot occur when
$\Pi=\Delta(C)$ for some linear completely regular code $C$.
\item[D.] We wonder whether it is true that if the quotient graph $\Gamma/\Pi$ in Theorem \ref{Tarithmeticquo} is isomorphic to a Hamming
graph or a Doob graph, then each cell in $\Pi$ can be expressed as
a cartesian product of completely regular codes with covering
radius at most two. We will show in Theorem \ref{Tproductrho=1} below that this holds when $\Pi$ is the coset partition of a linear completely regular code satisfying the eigenvalue
conditions of Theorem \ref{Tarithmeticquo}.
\end{description}
\end{remark}

In the special case when no two adjacent vertices are in the same cell of partition $\Pi$, we
can strengthen Theorem \ref{Tarithmeticquo} by looking at the images of the cliques in $\Gamma$.

\begin{prop}
\label{Pcliques}
Let $\Gamma$ be the Hamming graph $H(n,q)$. Let
$\Pi$ be any completely regular partition of $\Gamma$ where each
code $C$ in $\Pi$ has minimum distance $\delta(C)\geq 2$.
\begin{enumerate}
\item If $\Gamma/\Pi\cong H(m,q')$ then $q'\geq q$.
\item If $q\geq 4$ then $\Gamma/\Pi$ is not isomorphic to any Doob graph.
\item If $q\geq 3$ then $\Gamma/\Pi$ can not have the same
intersection array as any folded cube of diameter at least two
(including $\{6,5,4;1,2,6\}$).
\item Suppose further that $\delta(C)\ge 3$, and that  $\Gamma/\Pi$ has intersection array $\{6,5,4;1,2,6\}$.
Then $q=2$
and $\Gamma/\Pi$ is the folded $6$-cube.
\end{enumerate}
\end{prop}

\begin{proof}
Since each $C$ in $\Pi$ satisfies $\delta(C)\ge 2$, the vertices
in any clique in $\Gamma$ belong to pairwise distinct cells in
$\Pi$. So, in the quotient $\Gamma/\Pi$, every edge lies in a clique of size at least
$q$, the clique number of $H(n,q)$. This immediately implies
(a)--(c) as: $H(m,q')$ has clique number $q'$; any
Doob graph has maximal cliques of size three; a folded cube
other than $K_4$ has clique number two, as does any graph
with intersection array $\{6,5,4;1,2,6\}$.

For part (d), we clearly need only consider the case where $q=2$.
As $\delta(C)\geq 3$, we use \cite[Theorem~11.1.8]{B} to determine $n=\delta(C)=6$ from 
the intersection array of $\Gamma/\Pi$, and this means that $\Gamma/\Pi$ is the folded $6$-cube.
\end{proof}

As a side remark, we note that  every perfect $1$-code in $H(n,q)$, additive or not, 
gives us an equitable partition --- the partition into the code and its translates by the
various vectors of Hamming weight one ---  with complete quotient graph.  (This was a small 
oversight in \cite{B}, near the bottom of p355.)

\begin{center}
\begin{figure}
\includegraphics[angle=0,height=2.4in]{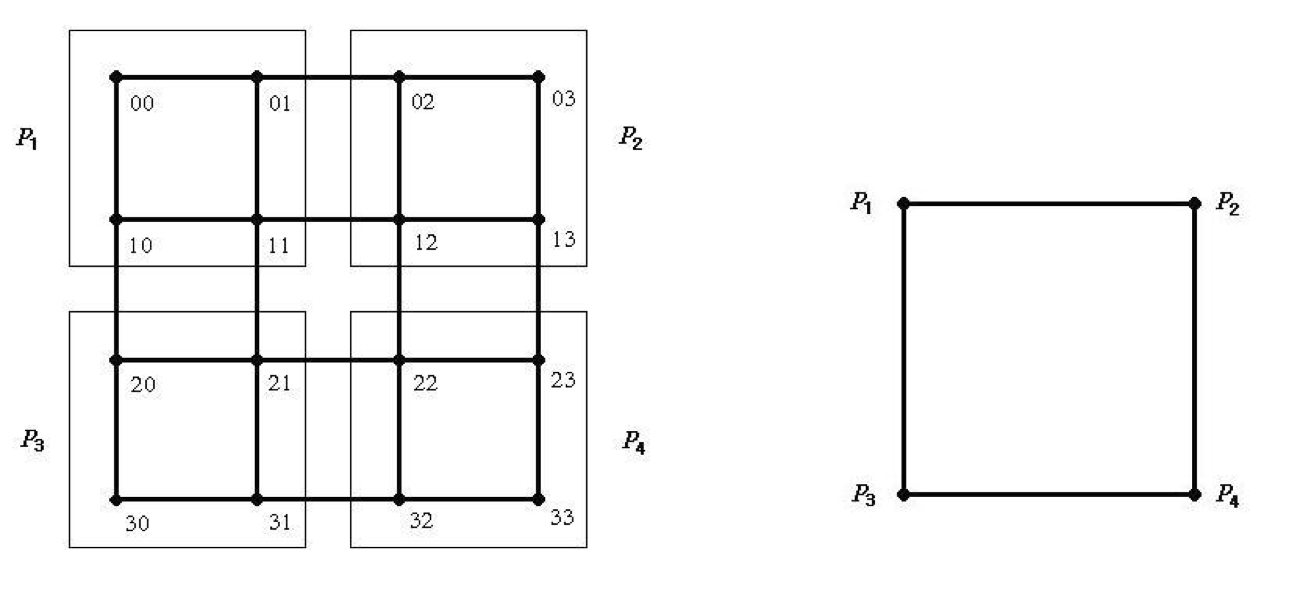}  
\caption{Hamming graph $H(2,4)$ and its quotient graph $H(2,2)$}
\end{figure}
\end{center}

\begin{Ex}
\label{EH24H22} Let $\Gamma$ be the Hamming graph $H(2,4)$ and let
$\Pi=\Pi(\Gamma)=\{P_1,P_2,P_3,P_4\}$, where
$P_1=\{00,01,10,11\}$, $P_2=\{02,03,12,13\}$,
$P_3=\{20,21,30,31\}$, and $P_4=\{22,23,32,33\}$. Then the
quotient graph $\Gamma/\Pi$ is the Hamming graph $H(2,2)$ as depicted in
Figure 1. Observe that each cell in partition $\Pi$ is
expressible as a cartesian product:
\begin{gather*}
	P_{1}=\{0,1\}\times\{0,1\}, \quad P_{2}=\{0,1\}\times\{2,3\}, \\
	P_{3}=\{2,3\}\times\{0,1\}, \quad P_{4}=\{2,3\}\times\{2,3\}.
\end{gather*}
Finally note that this example can be easily extended to give
$H(n,q)$ as a quotient of $H(n,sq)$, $n,s\geq 1$.
\end{Ex}

From Example \ref{EH24H22}, we see that the $n$-cube and the
folded $n$-cube can be quotients of $H(n,4)$. So Proposition
\ref{Pcliques}(a) does not hold when $\delta(C)=1$, even for additive
codes.

From now on, we will focus on linear completely regular codes over $GF(q)$.
We begin with a comment.
Recall that $C$ reduced means that it is impossible to express $C$ in the form $\mathcal{Q}\times C'$, where $\mathcal{Q}=GF(q)$ in this case.
Let $\ve_i$ denote the vertex of $H(n,q)$ with $i^{\mathrm{th}}$ position  $1$  and all other positions $0$.

\begin{lemma}
\label{Ldelta>=2}
Let $C$ be a
non-trivial linear completely regular code over $GF(q)$ of
length $n$.
Then $C$ is reduced if and only if $\delta(C)\geq 2$.
\end{lemma}

\begin{proof}
If $C$ contains a vector of weight one, say $\zeta \ve_{i}$, then by linearity $\ve_{i}$ is in $C$ and hence $C$ is not reduced.
The converse is also obvious. 
\end{proof}

We next show that, if $\Pi$ is the coset partition of
$H(n,q)$ with respect to a linear code with a quotient of Hamming
type, then each cell in $\Pi$ may be expressed as a cartesian
product of completely regular codes with covering radius one.

\begin{theorem}
\label{Tbreakdown}
Let $\Gamma$ be the Hamming graph $H(n,q)$. Let
$C$ be a linear completely regular code with minimum distance
$\delta(C)\geq 2$ in $\Gamma$ and let $\Delta(C)$ be the coset
partition of $\Gamma$ with respect to  $C$. 
\begin{enumerate}
\item Suppose $\Gamma/\Delta(C)\cong H(m,q')$.
Then $m$ divides $n$ and $C\cong \prod_{i=1}^{m}C^{(i)}$, where each $C^{(i)}$ is a linear $q$-ary
completely regular code with covering radius one and length $n/m$;
\item If $\Gamma/ \Delta(C)$ is isomorphic to a folded $m$-cube, $m\geq 4$, then $q=2$ and
\begin{equation*}
	C\cong \nullsp  {\begin{bmatrix}& H_1&|& \cdots&|& H_{m-1} &|& J &\end{bmatrix}},
\end{equation*}
where $H_i$ is the $(m-1)\times \gamma_1$ matrix with row $i$ all ones and all other rows all zero and $J$ is the $(m-1)\times \gamma_1$ all ones matrix, where $\gamma_1=\gamma_{1}(C)$.
\end{enumerate}
\end{theorem}

\begin{proof}
Note that $\ell_i:=\{\zeta\ve_{i}\,|\, \zeta\in GF(q)\}$ ($1\leq i\leq n$) are the singular lines (cliques of size $q$) in $\Gamma$ containing $\zero$.
Since $\delta(C)\geq 2$, the quotient map $\pi:V\Gamma\rightarrow V(\Gamma/\Delta(C))$ sending $\x\in V\Gamma$ to $\x+C$,
maps each singular line in $\Gamma$ into some singular 
line in $\Gamma/\Delta(C)$.
Hence we can consider the equivalence relation $\equiv$ on $\{1,\dots,n\}$ given by $i\equiv j$ if $\pi(\ell_i)$ and $\pi(\ell_j)$ belong to the same singular line in $\Gamma/\Delta(C)$.

First we prove (a).
Since $C$ is linear, we have either $\pi(\ell_i)=\pi(\ell_j)$ or $\pi(\ell_i)\cap\pi(\ell_j)=\{\zero\}$ for $1\leq i,j\leq n$.
Since there are precisely $(q'-1)\gamma_1$ neighbors of $\zero$ which are mapped to a fixed singular line in $\Gamma/\Delta(C)$ containing $\pi(\zero)=C$, each equivalence class has size $(q'-1)\gamma_1/(q-1)=n/m$, showing that $m$ divides $n$.

Let $R_1,\dots,R_m$ denote the equivalence classes of $\equiv$.
For $1\leq i\leq m$, let $D_i:=\{\x\in V\Gamma\,|\,\supp(\x)\subseteq R_i\}$.

\begin{claim}\label{claim 1}
For any $\vd\in D_{i}$, $\pi(\vd)$ belongs to the singular line of $\Gamma/\Delta(C)$ corresponding to $R_i$.
\end{claim}

\begin{proof}[Proof of Claim \ref{claim 1}]
Write $\vd=\sum_{j\in\supp(\vd)}\zeta_j\ve_j$.
If $\pi(\zeta_j\ve_j)=\pi(\zeta_h\ve_h)$ for some distinct $j,h\in\supp(\vd)$, then $\zeta_j\ve_j-\zeta_h\ve_h\in C$ and we may replace $\vd$ by $\tilde{\vd}:=\vd-(\zeta_j\ve_j-\zeta_h\ve_h)$.
Note that $\wt(\tilde{\vd})<\wt(\vd)$.
Hence we may assume from the first that $\pi(\zeta_j\ve_j)\ne\pi(\zeta_h\ve_h)$ for distinct $j,h\in\supp(\vd)$.
We now prove the result by induction on $\wt(\vd)$.
The result is obvious when $\wt(\vd)=1$, so assume $\wt(\vd)\geq 2$.
Take distinct $j,h\in\supp(\vd)$, and let $\vd':=\vd-\zeta_j\ve_j$ and $\vd'':=\vd-\zeta_h\ve_h$.
Then we have $\pi(\vd')\ne\pi(\vd'')$; for otherwise $\zeta_j\ve_j-\zeta_h\ve_h=\vd''-\vd'\in C$, so $\pi(\zeta_j\ve_j)=\pi(\zeta_h\ve_h)$, a contradiction.
Now $\pi(\vd)$ is adjacent to both $\pi(\vd')$ and $\pi(\vd'')$ which, by induction, belong to the above singular line of $\Gamma/\Delta(C)$, hence $\pi(\vd)$ must also belong to the same singular line.
\end{proof}

For $1\leq i\leq m$, let $\Sigma_i$ be the subgraph of $\Gamma$ induced on $D_i$, so $\Sigma_i\cong H(n/m,q)$.
Let $C^{(i)}:=D_{i}\cap C$.
Note that $C^{(i)}\ne D_i$ since $C$ is reduced.

\begin{claim}\label{claim 2}
$C^{(i)}$ is a linear $q$-ary completely regular code with covering radius one in $\Sigma_i$ with $U(C^{(i)})=\left(\begin{smallmatrix} 0 & \beta \\ \gamma_1 & \beta-\gamma_{1}\end{smallmatrix}\right)$ where $\beta:=n(q-1)/m$.
\end{claim}

\begin{proof}[Proof of Claim \ref{claim 2}]
Let $\vd\in D_i$ with $\vd\not\in C^{(i)}$.
By Claim \ref{claim 1}, $\pi(\vd)$ is adjacent to $\pi(\zero)$; i.e., $d(\vd,C)=1$.
Let $\x\in \Gamma(\vd)\cap C$, so $\x=\vd+\zeta\ve_j$ for some $j$ and non-zero $\zeta\in GF(q)$.
Then we have $\pi(-\zeta\ve_j)=\pi(\vd)$ and hence $\pi(\ell_j)$ must be contained in the singular line of $\Gamma/\Delta(C)$ corresponding to $R_i$.
This shows $j\in R_i$ and we have $\Gamma(\vd)\cap C=\Gamma(\vd)\cap C^{(i)}$.
This gives us the result.
\end{proof}

Since $C$ is additive, we have $\prod_{i=1}^m C^{(i)}\subseteq C$.
Moreover, by Proposition \ref{Pprodcrc}, $\prod_{i=1}^m C^{(i)}$ is a completely regular
code with the same intersection numbers as $C$ (i.e., the intersection numbers of $H(m,q')$
scaled by  $\gamma_1$). So $C=\prod_{i=1}^m C^{(i)}$. 

The proof of (b) is almost identical with that of (a) above.
First note that, as the folded $m$-cube has singular line size two, it follows that $q=2$ and $|R_i|=\gamma_1=n/m$ ($1\leq i\leq m$).
Moreover, each $C^{(i)}$ is a binary linear completely regular code with covering radius one, length $\gamma_1$ and $U(C^{(i)})=\left(\begin{smallmatrix} 0 & \gamma_1 \\ \gamma_1 & 0 \end{smallmatrix}\right)$; i.e., $C^{(i)}$ is the binary even weight code.
Let $C':=\prod_{i=1}^{m}C^{(i)}\subseteq C$.
Then we have $\dim C=n-m+1$, $\dim C'=m(\gamma_1-1)=n-m$, and $\Gamma/\Delta(C')\cong H(m,2)$.
It follows that $C/C'$ is the repetition code in $H(m,2)$ and hence $C=C'\cup(\x+C')$ for any word $\x$ with odd weight inside each equivalence class $R_i$.
\end{proof}

The next result rules out Doob graphs as quotients. 

\begin{prop}
\label{PnotDoob}
If $C$ is a linear $q$-ary completely regular
code of length $n$, then $H(n,q)/\Delta(C)$ is not isomorphic to
a Doob graph.
\end{prop}

\begin{proof} 
Since the Doob graph has $2^{2m}$ vertices for some $m$, we restrict to the
case where $q$ is a power of $2$.
Without loss of generality, we may assume $C$ is non-trivial and reduced, so $\delta(C)\geq 2$ by Lemma \ref{Ldelta>=2}.
Now if $H(n,q)/\Delta(C)$ is a Doob graph, then by Proposition
\ref{Pcliques}(b), we have $q\leq 3$, i.e., $q=2$.
Moreover, every vertex in $H(n,q)/\Delta(C)$ is contained in a maximal clique of size three.
Let $\{C,C_1,C_2\}$ be such a clique.
We may assume $\ve_1\in C_1$ and $\ve_2\in C_2$.
Then
some neighbor of $\ve_1$, say $\ve_1+\ve_i$, belongs to $C_2$. By linearity,
$\ve_1+\ve_2+\ve_i \in C$, so we have $i\ne 1,2$, and $\ve_i+C$ is adjacent in
$H(n,q)/\Delta(C)$ to all of $C,C_1,C_2$, contradicting our choice of a maximal clique.
\end{proof}

The following result was also shown by Borges et al. \cite{BRZ}. We give a proof for the convenience of the reader.

\begin{theorem}
\label{rho=1,2}
Let $C$ be a linear completely regular code over $GF(q)$ with
minimum distance $\delta(C)\geq 3$ in $H(n,q)$.
\begin{enumerate}
\item If $\Ev(C)=\{n(q-1),n(q-1)-qt\}$ for some $t\geq 1$, then
$C$ is a Hamming code. 
\item If $\Ev(C)=\{n(q-1),n(q-1)-qt,n(q-1)-2qt\}$ for some $t\geq 1$,
then $C=D\times D$ where $D$ is a Hamming code, or $q=2$ and $C$ is an extended Hamming code, or $q=2$, $n=5$ and $C$ is the repetition code. 
\end{enumerate}
\end{theorem}

\begin{proof}
(a) Since $C$ has covering radius one and minimum distance three,
it follows immediately that $C$ has to be perfect.

(b) Assume first that $t=1$.
Then, by a result of Meyerowitz \cite[Theorem 7]{BGKM}, $C$ takes the form $\{0\}^{\ell} \times \mathcal{Q}^{n-\ell}$ for some $0\le \ell \le n$ where $\mathcal{Q}=GF(q)$; so $C$ is either trivial or reduced, a contradiction.
Hence we have $t\ge 2$.

The coset graph $\Sigma:= H(n,q)/\Delta(C)$ is a connected strongly regular graph (i.e., a distance-regular graph with diameter two); consider the associated parameters $k,\lambda,\mu$ in the standard notation so that $\iota(\Sigma)=\{k,k-1-\lambda;1,\mu\}$. 
Note that $\mu>0$ since $\Sigma$ is connected.
Using \cite[Theorem 1.3.1]{B}, we obtain
\begin{align*}
	k&=n(q-1), \\
	\lambda&=(n(q-1)-qt)(n(q-1)-2qt+3), \\
	\mu&=n(q-1)+(n(q-1)-qt)(n(q-1)-2qt).
\end{align*}
Since $\mu\leq k$, we have $\frac{qt}{q-1}\leq n\leq\frac{2qt}{q-1}$.
This, combined with the constraint $0\leq\lambda\leq k-2$, gives us three possibilities:
(i) $n=\frac{qt}{q-1}$; (ii) $n=\frac{2qt-3}{q-1}$; (iii) $n=\frac{2qt-2}{q-1}$.

For (i), we have $k=qt$, $\lambda=0$, $\mu=qt$, and hence $\Sigma$ is the complete bipartite graph $K_{qt,qt}$.
Note that $2qt$ must divide $q^n$, so we have $q-1=1$, i.e., $q=2$, since $n$ is an integer.
Moreover, $C$ is an even weight code; for otherwise $\Sigma$ would have a cycle of odd length.
In particular, we have $\delta(C)\geq 4$.
By a result of Brouwer \cite{Brouwer}, it follows that a truncation of $C$ is a linear completely regular code with covering radius one, so it is a perfect code with minimum distance three.
Hence $C$ is an extended Hamming code.

For (ii), we have $0<\mu=6-qt$.
Since $t\geq 2$, we have $q=t=2$.
So $n=5$, $k=5$, $\lambda=0$, $\mu=2$, and hence $\Sigma$ has the intersection array $\{5,4;1,2\}$ of the folded $5$-cube.
By \cite[Theorem 9.2.7]{B}, $\Sigma$ is isomorphic to the folded $5$-cube.
By Theorem \ref{Tbreakdown}(b) and since $\gamma_1=1$, it follows that $C$ is the repetition code with length five.

For (iii), we have $k=2qt-2$, $\lambda=qt-2$, $\mu=2$, so $\Sigma$ has the intersection array $\{2qt-2,qt-1;1,2\}$ of the Hamming graph $H(2,qt-1)$.
By \cite{E}, $\Sigma$ is isomorphic to $H(2,qt-1)$.
(Note that $qt-1\ne 4$ as $q\geq 2$ and $t\geq 2$.)
It follows from Theorem \ref{Tbreakdown}(a) that $C$ is the cartesian product of two copies of the Hamming code with length $n/2$.
\end{proof}

Let $C$ be a non-trivial, reduced linear completely regular code over $GF(q)$
in $\Gamma=H(n,q)$ with intersection numbers $\gamma_i,\alpha_i,\beta_i$ and parity check matrix $H$.
Note that the columns of $H$ are all non-zero since $C$ is reduced.
Consider the equivalence relation $\simeq$ defined on $\{1,\dots,n\}$ by\footnote{With the notation in the proof of Theorem \ref{Tbreakdown}, $i\simeq j$ if and only if $\pi(\ell_i)=\pi(\ell_j)$.} $i \simeq j$ if the $i^{\mathrm{th}}$ and $j^{\mathrm{th}}$ columns of $H$ are linearly dependent, or equivalently, $i=j$ or $\{i,j\}=\supp(\x)$ for some $\x\in C$.
Note that such $\x$ is unique up to scalar multiplication; for otherwise we would have $\ve_i,\ve_j\in C$ and hence $C$ is non-reduced, a contradiction.
One also easily checks that all equivalence classes have size $\gamma_1$. 

\begin{lemma}
\label{LlittleDcrc}
With the above notation, let $P\subseteq  \{1,\dots,n\}$ be a complete set of representatives for the equivalence classes of $\simeq$.
Let $D:=\{\x\in V\Gamma\,|\,\supp(\x)\subseteq P\}$ and let $\Sigma\cong H(n/\gamma_1,q)$ be the subgraph of $\Gamma$ induced on $D$.
Then $C_P:=C\cap D$ is a linear completely regular code in $\Sigma$
with $\delta(C_P)\geq 3$, $\rho(C_P)=\rho(C)$, and intersection numbers $\gamma_i/\gamma_1,\alpha_i/\gamma_1,\beta_i/\gamma_1$.
Moreover, $\Gamma/\Delta(C)$ and $\Sigma/\Delta(C_P)$ are isomorphic.
\end{lemma}

\begin{proof}
Let $\Pi(C)=\{C_0=C,C_1,\dots,C_{\rho}\}$ be the distance partition of $V\Gamma$ with respect to $C$, where $\rho:=\rho(C)$.
We claim that $\{C_0\cap D=C_P,C_1\cap D,\dots,C_{\rho}\cap D\}$ is the distance partition of $V\Sigma$ with respect to $C_P$.
Indeed, let $\x\in C_i\cap D$, and suppose $\x+\zeta_i\ve_i\in C_{i-1}$ for some non-zero $\zeta_i\in GF(q)$.
Then for any $j$ with $i\simeq j$, there is a unique non-zero $\zeta_j\in GF(q)$ such that $\zeta_i\ve_i-\zeta_j\ve_j\in C$, so $\x+\zeta_j\ve_j=\x+\zeta_i\ve_j-(\zeta_i\ve_i-\zeta_j\ve_j)\in C_{i-1}$.
Hence it follows that $\x$ has precisely $\gamma_i/\gamma_1$ neighbors in $C_{i-1}\cap D$.
By the same reasoning, $\x$ has $\alpha_i/\gamma_1$ neighbors in $C_i\cap D$ and $\beta_i/\gamma_1$ neighbors in $C_{i+1}\cap D$.
In particular, $C_P$ is a completely regular code in $\Sigma$ with $\rho(C_P)=\rho$.
Since any word in $V\Gamma$ is congruent modulo $C$ to some word in $D$, it also follows that $\Gamma/\Delta(C)$ and $\Sigma/\Delta(C_P)$ are (canonically) isomorphic.
Finally, since $C_P$ can not contain words with weight two, we have $\delta(C_P)\geq 3$.
\end{proof}

\begin{theorem}
\label{Tproductrho=1} Let $C$ be a non-trivial, reduced, linear
completely regular code over $GF(q)$ in $H(n,q)$ with intersection number $\gamma_1=\gamma_1(C)$.
Suppose that $\Ev(C)=\{n(q-1),n(q-1)-qt,\dots,n(q-1)-q\rho t \}$ for some $t$.
Then one of the following holds:
\begin{enumerate}
\item $q=2$ and
\begin{equation*}
	C\cong \nullsp \underbrace{\begin{bmatrix}&M &|&\cdots&|&M &\end{bmatrix}}_{\gamma_1 \, \text{copies}},
\end{equation*}
where $M=[ I \,|\, \one]$  is a parity check matrix for a binary repetition code, and the quotient is a folded cube;
\item $\rho=1$ and
\begin{equation*}
	C\cong \nullsp \underbrace{\begin{bmatrix}&H &|&\cdots&|&H&\end{bmatrix}}_{\gamma_1\, \text{copies}},
\end{equation*}
where $H$ is a parity check matrix for some Hamming code, and the quotient is a complete graph;
\item $\rho=2$, $q=2$ and
\begin{equation*}
	C\cong \nullsp \underbrace{\begin{bmatrix}&E&|&\cdots&|&E&\end{bmatrix}}_{\gamma_1\, \text{copies}},
\end{equation*}
 where $E$ is a parity check matrix for a fixed extended Hamming code, and the quotient is a regular complete bipartite graph;
\item $\rho\geq 2$ and
\begin{equation*}
	C\cong\underbrace{C_{1}\times\dots\times C_{1}}_{\rho\, \text{copies}},
\end{equation*}
where $C_{1}$ is a completely regular code with covering radius one, and the quotient is a Hamming graph.
\end{enumerate}
\end{theorem}

\begin{proof}
The reader can easily check that examples (a)--(d)  are all
completely regular. (See also Bier \cite{Bier}.)

Recall the code $C_P$ defined in Lemma \ref{LlittleDcrc}, which is a linear completely regular code in $\Sigma\cong H(n/\gamma_1,q)$ with $\delta(C_P)\geq 3$ and still satisfies our arithmetic hypothesis.
Note that if $\delta(C)\geq 3$ then we have $\gamma_1=1$ so $C=C_P$.
Note also that if $K$ is a parity check matrix for $C_P$, then we have $C \cong \nullsp \left[ \begin{array}{c|c|c} K & \cdots & K \end{array} \right]$ ($\gamma_1$ copies).
Hence we may restrict to the case where $\delta(C)\geq 3$.

The result follows from Theorem \ref{rho=1,2} when $\rho\in\{1,2\}$, so we assume $\rho\geq 3$.
Then Theorem \ref{Tarithmeticquo} tells us the possible quotients, i.e., a folded cube, a Hamming graph, a Doob graph, or the incidence graph of a symmetric $2$-$(16,6,2)$ design.
If the coset graph is a Hamming graph, then (d) applies by Theorem \ref{Tbreakdown}(a).
If it is a folded $m$-cube, then Theorem \ref{Tbreakdown}(b) gives us a parity check matrix for $C$ and one easily sees that rearranging columns of this matrix, we obtain the parity check matrix given in (a).
By Proposition \ref{PnotDoob}, a Doob graph cannot arise as coset graph.
Finally, if the coset graph has intersection array $\{6,5,4;1,2,6\}$, then it is indeed the folded $6$-cube by Proposition \ref{Pcliques}(d) and we are done.
\end{proof}

We remark that the above result includes the following generalization
of the result of Bier \cite{Bier} concerning coset graphs which
are isomorphic to Hamming graphs (cf.~\cite[p354]{B}):

\begin{cor}
\label{CHammingquo} Let $C$ be a reduced linear completely regular code in
$H(n,q)$ with intersection number $\gamma_1=\gamma_1(C)$, whose coset graph is a Hamming graph $H(m,q')$.
Then one of the following holds:
\begin{enumerate}
\item $\rho=1$ and
\begin{equation*}
	C\cong \nullsp {\underbrace{\begin{bmatrix}&H &|&\cdots&|&H&\end{bmatrix}}_{\gamma_1\, \text{copies}}},
\end{equation*}
where $H$ is a parity check matrix for some Hamming code;
\item $\rho\geq 2$ and
\begin{equation*}
	C\cong\underbrace{C_{1}\times\cdots\times C_{1}}_{\rho\, \text{copies}},
\end{equation*}
where $C_{1}$ is a completely regular code with covering radius one.
\end{enumerate}
\end{cor}

\begin{proof}
The eigenvalues of $H(m,q')$ are in arithmetic progression, and hence so are the eigenvalues of $C$ in view of Proposition \ref{Pbcnquo}.
Now the result follows from Theorem \ref{Tproductrho=1}.
\end{proof}

\section*{Acknowledgments}
This work is part of the Ph.D. thesis of WSL under supervision of JHK.
Part of this work was completed while JHK was with the Pohang University of Science and Technology (POSTECH) and  WJM was visiting him.  WJM
wishes to thank the Department of Mathematics at POSTECH for their
hospitality and Com$^2$MaC for financial support. JHK is partially supported by the `100 talents' program of the Chinese government.
 WJM wishes to thank the US
National Security Agency for financial support under grant numbers
H98230-07-1-0025 and H98230-12-1-0243. 
HT is supported in part by JSPS KAKENHI Grant No.~25400034.
The authors thank the editor for patiently guiding the manuscript through the review process.


\newpage 

\section*{\large ADDENDUM: Correcting an error in Lemma 6.3 of \cite{KLM}}

We now point out an error in the statement of Lemma 6.3 in our companion 
paper \cite{KLM}; a key hypothesis is missing from the statement of the lemma.

We first recall several definitions from \cite{KLM}.

Let $\Gamma$ be a distance-regular graph.
We say that a subspace $W$ of $V:=\mathbb{C}^{V\Gamma}$ is {\em Schur closed} if $\u \circ \vv \in W$ whenever $\u$ and $\vv$ themselves
belong to $W$. The {\em Schur closure} of a subspace (or set) is the smallest Schur closed subspace containing it.
A subspace is easily seen to be Schur closed if and only if it contains some basis of pairwise orthogonal 01-vectors.

Let $E_0,E_1,\dots, E_D$ be any ordering of the
primitive idempotents of $\Gamma$. Let $V_j\subseteq V$ denote the eigenspace corresponding to $E_j$.
For a completely regular code $C \subseteq V\Gamma$ with $\Ev(C)=\{\theta_{i_0},\dots, \theta_{i_\rho}\}$, 
consider the characteristic vector $\x=\x_C\in V$ of $C$ and the {\em outer distribution module} $\bma\x$ which has 
$\left\{ E_{i_0} \x , \dots, E_{i_\rho} \x  \right\}$ as basis. It is well known that $\bma \x$ is Schur closed whenever
$C$ is completely regular.  We say \cite[Definition~4.1]{KLM} $C$ is {\em Leonard} 
(with respect to this ordering) if 
\begin{equation*}
\u^{(\ell)}:={\underbrace {\u\circ\dots\circ \u}_{\ell \ \text{times}}}\in
(V_{i_{0}}+ V_{i_1} +\dots+ V_{i_{\ell}}) \setminus (V_{i_0} + V_{i_1} +\dots+V_{i_{\ell - 1}}),
\end{equation*}
for $1\le \ell \le \rho$, where  $\u:= E_{i_1} \x$.
On the other hand, when $E_0,E_1,\dots, E_D$ is a $Q$-polynomial ordering, we say that completely regular code $C$ is
{\em harmonic} (with respect to this ordering) if $\Ev( C ) = \{ \theta_{ti} \,|\, 0\le i \le \rho \}$ for some $t$  \cite[Definition~6.1]{KLM}.

We now give the corrected statement of Lemma 6.3 in \cite{KLM} together with a full 
proof\footnote{We found a typographical error in the proof
of Lemma 6.3 in \cite{KLM}: where we write ``$\le t$'' and ``$\le 1$'' in the displayed equations, it should instead read ``$=t$'' and ``$=1$'', respectively.}.

\setcounter{section}{6}
\setcounter{lemma}{2}

\begin{lemma}
\label{Lemma 6.3}
Let $\Gamma$ be a distance-regular graph which is $Q$-polynomial with
respect to the ordering $E_{0},E_{1},\dots,E_{D}$  of its primitive idempotents. Let $C$ be a completely 
regular code which is harmonic with respect to this ordering with $\Ev(C)=\{ \theta_{ti} \,|\, 0\le i \le \rho\}$. 
If $E_t \x$ has $\rho+1$ distinct entries, then $C$ is a Leonard completely regular code. In particular, 
whenever the $Q$-polynomial ordering satisfies $\theta_0 > \theta_1 > \dots > \theta_D$,
every completely regular code which is harmonic with respect to this ordering is also Leonard 
with respect to this ordering.
\end{lemma}

\begin{proof}
The outer distribution module $\bma \x$ has $\left\{ E_{ti} \x \,|\, 0\le i \le \rho\right\}$ as a basis.
Since  $\bma\x$ is Schur closed, 
there exist numbers $\omega_{\ell,j}$ such that
$E_{t}\x \circ E_{tj}\x =\sum_{\ell=0}^{\rho}\omega_{\ell,j}E_{t \ell }\x$. 
Using a fundamental result of Cameron, Goethals and Seidel 
 (see, e.g., \cite[Theorem~2.2]{KLM}), we see that, since $\Gamma$ is $Q$-polynomial, 
 we must have $\omega_{\ell,j}  = 0$ whenever  $|\ell-j|>1$. 
 
 Now we employ the assumption (missing in the statement of Lemma~6.3 in \cite{KLM}) that $E_t \x$ has $\rho+1$
 distinct entries to show that $\omega_{\ell,j} \neq 0$ when $\ell=j+1$ for $0 \le  j<  \rho$. Under this assumption,
 the Schur closure of $\{ E_t \x\}$ has dimension $\rho+1$. If there were some $j<\rho$ with $\omega_{\ell,j}=0$ for all
 $\ell > j$, this would give us a proper Schur closed $\bma$-submodule $\left\{ E_{ti} \x \,|\, 0\le i \le j \right\}$ which is impossible.
Now let $\u := E_t \x$. Then we use $\omega_{\ell,j} =0$ for $\ell>j+1$ to see that 
 \begin{equation*}
 	\u^{(j+1)} \in V_{0}+ V_{t} +\dots+ V_{t(j +1) }
\end{equation*}
while 
\begin{equation*}
	\u^{(j+1)} \not\in V_{0}+ V_{t} +\dots+ V_{t j }
\end{equation*}
follows since $\omega_{j+1,j} \neq 0$ and, by induction, $E_{t j}  \u^{(j)} \neq 0$.
Thus the code $C$ is Leonard with respect to this ordering.
 
The most common $Q$-polynomial orderings known have $\theta_0>\theta_1> \dots > \theta_D$.
In this case $\theta_t$ is the second largest eigenvalue of the tridiagonal matrix $U=U(C)$ and standard results about Sturm sequences (see, e.g., \cite{Martin-crdt1,Martin-mindist}) guarantee that $E_t\x$ has $\rho+1$ distinct entries in all such cases.
\end{proof}
 
Of particular relevance here is the observation that the lemma applies to all arithmetic completely regular codes in the Hamming graphs.
For if $C$ is arithmetic, then it is harmonic, and thus Leonard, with respect to the natural $Q$-polynomial ordering for $H(n,q)$.
But what of the second $Q$-polynomial ordering in the case $q=2$?

In \cite{Dickie}, Dickie showed that the only distance-regular graphs with diameter $D\ge 5$ which admit more than one $Q$-polynomial
ordering are dual polar graph $[{}^2A_{2D-1}(r)]$ ($r\ge 2$, a prime power), the $D$-cube $H(D,2)$ (when $D$ is even), the halved $(2D+1)$-cube, the folded $(2D+1)$-cube, and the regular $n$-gons $n=2D,2D+1$. In all these cases, one must be careful in applying Lemma 6.3 in \cite{KLM}
when the $Q$-polynomial ordering does not correspond to the natural ordering of graph eigenvalues from largest to smallest.

We end with an example which underscores the 
necessity of the condition that $E_t\x$ have $\rho+1$ distinct entries.

The trivial completely regular codes $C=\{0,1 \}^{n-\rho} \times \{0\}^\rho$ in the $D$-cube $H(n,2)$ have 
$\Ev(C)=\{ n, n-2, n-4, \dots, n-2\rho\}$ so, with respect to the standard ordering, 
these codes are arithmetic (with $t=1$), hence both  harmonic and Leonard. 
But when $n$ is even, the $n$-cube admits a second
$Q$-polynomial structure, with eigenvalue ordering
\begin{equation*}
	\theta_j = (-1)^j (n-2j) \quad (0\le j \le n).
\end{equation*}
With respect to this ordering, $E_1 \x$ has only  $1 + \lfloor \rho/2 \rfloor$ distinct entries and  covering radius $\rho$. Under this ordering,  
\begin{equation*}
	\Ev(C)=\{ \theta_0,\theta_2,\dots, \theta_{2\lfloor \rho/2 \rfloor} , \ 
	\theta_{ n+1-2\lceil \rho/2 \rceil}, \dots,  \theta_{n-3}, \theta_{n-1} \}
\end{equation*}
 so that, in the particular case $\rho=n-1$, $C$ is arithmetic. Varying $n$, we find an 
 infinite family of examples (codes of size two) which are harmonic but not Leonard.

\end{document}